\documentclass{amsart}
\usepackage{amsmath,amsthm,amsfonts,amssymb,amsxtra,mathrsfs}
\usepackage{enumerate}
\usepackage[backref=page]{hyperref}
\hypersetup{colorlinks=true,citecolor=blue,linkcolor=blue,urlcolor=blue,pdfstartview=FitH,
pdfauthor=Tobias Berger and Gergely Harcos,pdftitle=l-adic
representations associated to modular forms over imaginary quadratic
fields}

\newtheorem{thm}{Theorem}[section]

\newtheorem{lem}[thm]{Lemma}

\theoremstyle{definition}
\newtheorem{defn}[thm]{Definition}
\newtheorem*{acks}{Acknowledgements}
\theoremstyle{remark}
\newtheorem{rem}[thm]{Remark}
\numberwithin{equation}{section}
\newcommand{\ov}[1]{\overline{#1}}

\newcommand{\C}{\mathbb{C}}
\newcommand{\Q}{\mathbb{Q}}
\newcommand{\Z}{\mathbb{Z}}
\newcommand{\A}{\mathbb{A}}
\newcommand{\TT}{\ov{\mathbb{T}}}
\newcommand{\HH}{\mathscr{H}}
\newcommand{\MM}{\mathcal{M}}
\newcommand{\OO}{\mathcal{O}}
\newcommand{\eps}{\varepsilon}

\DeclareMathOperator{\Ind}{Ind} \DeclareMathOperator{\disc}{disc}
\DeclareMathOperator{\Gal}{Gal} \DeclareMathOperator{\Frob}{Frob}
\DeclareMathOperator{\Hom}{Hom} \DeclareMathOperator{\End}{End}
\DeclareMathOperator{\GL}{GL} \DeclareMathOperator{\GSp}{GSp}
\DeclareMathOperator{\GO}{GO} \DeclareMathOperator{\diag}{diag}

\begin{document}

\title[$\ell$-adic representations over imaginary quadratic fields]
{$\ell$-adic representations associated to modular forms over
imaginary quadratic fields}

\author{Tobias Berger}
\address{Department of Pure Mathematics and Mathematical Statistics,
Centre for Mathematical Sciences, Cambridge, CB3 0WB, United
Kingdom} \email{t.berger@dpmms.cam.ac.uk}

\author{Gergely Harcos}
\address{Alfr\'ed R\'enyi Institute of Mathematics, Hungarian Academy of Sciences, POB 127, Budapest H-1364, Hungary} \email{gharcos@renyi.hu}

\begin{abstract}
Let $\pi$ be a regular algebraic cuspidal automorphic representation
of $\GL_2$ over an imaginary quadratic number field $K$, and let
$\ell$ be a prime number. Assuming the central character of $\pi$ is
invariant under the non-trivial automorphism of $K$, it is shown
that there is a continuous irreducible $\ell$-adic representation
$\rho$ of $\Gal(\ov K/K)$ such that $L(s,\rho_v)=L(s,\pi_v)$
whenever $v$ is a prime of $K$ outside an explicit finite set.
\end{abstract}

\keywords{Galois representations, automorphic representations,
Siegel modular forms, Langlands conjectures}


\subjclass[2000]{Primary 11F80; Secondary 11F46, 11F70, 11R39}

\maketitle
\tableofcontents


\section{Introduction}

Let $K$ be an imaginary quadratic field with non-trivial
automorphism $c$, and let $\pi$ be a cuspidal automorphic
representation of $\GL_2(\A_K)$ with unitary central character
$\omega$. If $\pi_\infty$ has Langlands parameter
$W_\C=\C^\times\to\GL_2(\C)$ given by $z\mapsto\diag(z^{1-k},{\ov
z}^{1-k})$ for some integer $k\geq 2$ (that is, in the sense of
Clozel \cite{C}, $\pi$ is any regular algebraic cuspidal automorphic
representation up to twist), then by the Langlands philosophy $\pi$
should give rise (for any prime number $\ell$) to a continuous
irreducible $\ell$-adic representation $\rho$ of the Galois group
$\Gal(\ov K/K)$ such that the associated $L$-functions agree. In
other words, at each prime $v$ of $K$ the Frobenius polynomial of
$\rho$ at $v$ agrees with the Hecke polynomial of $\pi$ at $v$.
Under the assumption that $\omega=\omega^c$ it is possible to relate
$\pi$ to holomorphic Siegel modular forms via theta lifts and deduce
(using $\ell$-adic cohomology on Siegel threefolds) some weak
version of this predicted correspondence. In fact Taylor \cite{T}
managed to obtain the above equality of Frobenius and Hecke
polynomials for all $v$ outside a zero density set of places, but he
had to make some additional technical assumptions. It is our aim
here to describe how the results of Friedberg and Hoffstein
\cite{FH} on the non-vanishing of certain central $L$-values and
those of Laumon \cite{L1,L3} and Weissauer \cite{We} on associating
Galois representations to Siegel modular forms enable one to remove
these technical assumptions and conclude the statement for all $v$
outside an explicit finite set. Precisely, we shall prove the
following

\begin{thm}\label{thm1} Assume that $\omega=\omega^c$.
Let $S$ denote the set of places in $K$ which divide $\ell$ or where
$K/\Q$ or $\pi$ or $\pi^c$ is ramified. There exists a continuous
irreducible representation $\rho:\Gal(\ov K/K)\to \GL_2(\ov
\Q_\ell)$ such that if $v$ is a prime of $K$ outside $S$ then $\rho$
is unramified at $v$ and the characteristic polynomial of
$\rho(\Frob_v)$ agrees with the Hecke polynomial of $\pi$ at $v$,
that is\footnote{Here and later we shift-normalize all $L$-functions
such that $s=1/2$ is their center.}, $L(s,\rho_v)=L(s,\pi_v)$.
\end{thm}

\begin{rem} This theorem strengthens Theorem~A of \cite{T}.
The assumption $\omega=\omega^c$ is inherent to Taylor's method; one
would hope to remove this condition by other methods.\end{rem}

\begin{rem} By \cite[Lemma~6]{T2} the Galois representation takes values
in ${\rm GL}_2(E)$ for a finite extension $E/\Q_{\ell}$ (see
\cite[p.~635]{T}  for an explicit description).\end{rem}

\begin{rem}\label{known} There are two cases in which Theorem~\ref{thm1} is known to
hold (cf. \cite[Lemmata~1 and 2]{T}):
\begin{enumerate}
\item[1.]  $\pi\otimes\delta\cong\pi$ for some nontrivial quadratic character $\delta$ of $K$.
In this case if $L/K$ denotes the quadratic extension corresponding
to $\delta$ then there is an algebraic id\`ele class character
$\psi$ of $L$ unequal to its conjugate under the non-trivial element
of $\Gal(L/K)$ such that $\pi$ is the automorphic induction of
$\psi$ to $K$; the conclusion of Theorem~\ref{thm1} follows by work
of Serre \cite{S}.
\item[2.] $\pi\otimes\nu\cong(\pi\otimes\nu)^c$ for some finite order
character $\nu$ of $K$. In this case a twist of $\pi$ is a base
change from $\Q$; the conclusion of Theorem~\ref{thm1} follows from
work of Deligne \cite{D}.
\end{enumerate}
\end{rem}

The proof of Theorem~\ref{thm1} can be briefly outlined as follows.
The initial strategy is that of Taylor \cite{T}. We assume we are
not in a case covered by Remark~\ref{known}. Using the deep results
of \cite{HST} and \cite{FH} we construct a nonzero theta lift on
$\GSp_4(\A_\Q)$ of the twist $\pi\otimes\mu$ for a dense set of
quadratic id\`ele class characters $\mu$ of $K$ (see
Definition~\ref{defn1} below). The irreducible constituent $\Pi^\mu$
of such a lift is generated by a vector-valued holomorphic
semi-regular cusp form on the Siegel three-space. Using Hasse
invariant forms and the theory of pseudo-representations developed
by Wiles \cite{W} and Taylor \cite{T2,T3}, Taylor had shown that one
can associate a $4$-dimensional representation to $\Pi^\mu$ if one
could associate $4$-dimensional Galois representations to regular
holomorphic Siegel cusp forms. This is now possible by work of
Laumon \cite{L1,L3} and Weissauer \cite{We}. We obtain therefore,
for each $\mu$ in some dense set, a $4$-dimensional representation
of $\Gal(\ov \Q/\Q)$ with the same partial $L$-function as the one
associated to $\Pi^\mu$, and we prove that it is induced from some
$2$-dimensional representation $\rho^\mu$ of $\Gal(\ov K/K)$. By
exploring global compatibility relations among the various
$\rho^\mu$ we show that they can be replaced by quadratic
twists\footnote{Here and later we identify finite order id\`ele
class characters with continuous Galois characters.}
$\rho\otimes\mu$ of a single $2$-dimensional representation $\rho$
of $\Gal(\ov K/K)$, and we verify that this $\rho$ has the required
property of Theorem~\ref{thm1}.

\begin{rem} Since the construction of the $4$-dimensional Galois representation
associated to $\Pi^\mu$ involves an $\ell$-adic limit process one
loses information about the geometricity of the Galois
representation. For $\ell$ split in $K/\Q$, Urban \cite[Corollaire
2]{U} has proved, however, that if $\pi$ is ordinary then the Galois
representation $\rho$ of Theorem~\ref{thm1} is ordinary at $v \mid
\ell$.
\end{rem}

\section{Theta lifts}
Theorem~1 and Proposition~5 of \cite{HST} show how to construct
non-zero theta lifts on ${\rm GSp}_4$ of many quadratic twists of
$\pi$, conditional on ``Conjecture/Theorem~1" \cite[p.~403]{HST}.
The analytic non-vanishing result of \cite{FH} implies
``Conjecture/Theorem~1". For completeness, we decided to summarize
how \cite{FH} implies a strengthening of Theorem~1 of \cite{HST}.

By assumption the central character of $\pi$ factors through the
norm map as $\omega=\tilde{\omega}\circ N_{K/\Q}$, where
$\tilde{\omega}$ is a character of $\Q$. The ratio of the two
characters $\tilde\omega$ satisfying this equation is the quadratic
character corresponding to $K/\Q$, hence one of them is odd and the
other one is even. We shall consider the character $\tilde\omega$
with $\tilde\omega_\infty(-1)=(-1)^k$. By Proposition 1 of
\cite{HST} the pair $(\pi,\tilde{\omega})$ defines a cuspidal
automorphic representation of $\GO^\circ(\A_\Q)$, where $\GO$ is the
group of orthogonal similitudes of a certain quadratic space $W_K$
of sign $(3,1)$ over $\Q$. \cite{HST} also introduces a signature
$\delta=(\delta_v)$, a map from the places of $\Q$ to $\{\pm 1\}$
which is 1 at all but finitely many places such that $\delta_v=1$
whenever $\pi_v \not \cong \pi_v^c$ (here we view $\pi$ as a
representation of $R_{K/\Q}\GL_2$, the group obtained from $\GL_2$
by restriction of scalars). By Proposition 2 of \cite{HST} the
triple $\hat\pi:=(\pi,\tilde{\omega},\delta)$ modulo the action of
$\{1,c\}$ can be identified with a cuspidal automorphic
representation of $\GO(\A_\Q)$ which in turn has a theta lift
$\Theta(\hat\pi)$ to $\GSp_4(\A_\Q)$. By Proposition 3 of \cite{HST}
the lift $\Theta(\hat\pi)$ is contained in the space of cusp forms,
and if $\Pi$ denotes an irreducible constituent of $\Theta(\hat\pi)$
then $\Pi_\infty$ is a holomorphic limit of discrete series
representation of weight $(k,2)$ whenever $\delta_\infty=-1$, while
$\Pi_v$ is an unramified irreducible principal series representation
with $L(s,\Pi_v)=L(s,\pi_v)$ whenever $\delta_v=1$ and $v$ is a
rational prime which does not lie under a prime in $S$. In addition,
$\Pi$ is nonzero assuming there is a character $\varphi$ of $K$
restricting to $\tilde{\omega}$ on $\Q$ and satisfying the following
two properties:
\[\pi_v\cong \pi_v^c \quad\Longrightarrow\quad
\delta_v=\tilde{\omega}_v(-1)\eps(\pi_v\otimes\varphi_v^{-1},1/2);\]
\[L(\pi\otimes\varphi^{-1},1/2)\neq 0.\]

Using these as preliminaries we can deduce from the non-vanishing
result \cite{FH} that $\pi\otimes\mu$ gives rise to suitable
$\Pi^\mu$ on $\GSp_4(\A_\Q)$ for a dense set of quadratic characters
$\mu$ of $K$.

\begin{defn}\label{defn1} A set $\MM$ of quadratic characters of $K$
is {\it dense} if it has the following property. If $\tilde\mu$ is a
quadratic character of $K$ and $M$ is a finite set of rational
primes then there is a character $\mu\in\MM$ such that
$\mu_v=\tilde\mu_v$ for all $v\in M$.
\end{defn}

\begin{defn}\label{defn2} For a cuspidal automorphic representation $\tilde\pi$ of $K$
let $S_K(\tilde\pi)$ denote the set of places in $K$ which divide
$\ell$ or where $K/\Q$ or $\tilde\pi$ or $\tilde\pi^c$ is ramified,
and let $S_\Q(\tilde\pi)$ denote the set of rational places which
lie under some place in $S_K(\tilde\pi)$.
\end{defn}

\begin{thm}\label{nthm1} There exists a dense set $\MM$ of quadratic
characters of $K$ with the following property. For each $\mu\in\MM$
there is a signature $\delta$ such that the representation
$(\pi\otimes\mu,\tilde\omega,\delta)$ of $\,\GO(\A_\Q)$ gives rise
to a cuspidal automorphic representation $\Pi^\mu$ of
$\,\GSp_4(\A_\Q)$ satisfying:
\begin{itemize}
\item $\Pi_\infty^\mu$ is a holomorphic limit of discrete series representation of
weight $(k,2)$;
\item if $v$ is a rational prime outside $S_\Q(\pi\otimes\mu)$ then
$\Pi_v^\mu$ is an unramified irreducible principal series
representation with $L(s,\Pi_v^\mu)=L(s,(\pi\otimes\mu)_v)$.
\end{itemize}
\end{thm}

To prove this let $\tilde\mu$ be a quadratic character of $K$ and
let $M$ be a finite set of rational places. In the light of the
above discussion (i.e.\ by Propositions~1--3 of \cite{HST}) it
suffices to show that for $\tilde\pi:=\pi\otimes\tilde\mu$ there
exist a quadratic character $\eta$ of $K$ with $\eta_v=1$ for all
$v\in M$, a signature $\delta$ with $\delta_\infty=-1$ and
$\delta_v=1$ for any rational prime $v\notin
S_\Q(\tilde\pi\otimes\eta)$, and a character $\varphi$ of $K$ with
$\varphi|_\Q=\tilde{\omega}$, satisfying the additional properties
\begin{equation}\label{eq1}
\delta_v=\begin{cases}
\tilde{\omega}_v(-1)\eps(\tilde\pi_v\otimes\eta_v\varphi_v^{-1},1/2)
&\text{if $\tilde\pi_v\otimes\eta_v\cong (\tilde\pi_v\otimes\eta_v)^c$},\\
1&\text{if $\tilde\pi_v\otimes\eta_v\not\cong
(\tilde\pi_v\otimes\eta_v)^c$};
\end{cases}
\end{equation}
\begin{equation}\label{eq2}
L(\tilde\pi\otimes\eta\varphi^{-1},1/2)\neq 0.
\end{equation}
The proofs of Lemma~13 and Proposition 5 in \cite{HST} provide us
with $\eta$ and $\varphi$ satisfying
\begin{equation}\label{neq1}
\eps(\pi_\infty\otimes\varphi_\infty^{-1},1/2)
=-\tilde{\omega}_\infty(-1)
\end{equation}
and
\begin{equation*}
\eps(\tilde\pi\otimes\eta\varphi^{-1},1/2)=1.
\end{equation*}
Here we used our initial assumptions for $\pi$ and $\tilde\omega$.
Theorem~A and the first part of Theorem~B in \cite{FH} show that
$\eta$ can be replaced by another quadratic character satisfying
\eqref{eq2}. Now we define $\delta$ according to \eqref{eq1}. By
\cite[Lemma~14]{HST}, $\delta_v=\pm 1$ for all rational primes, and
$\delta_v=1$ for all rational primes $v\notin
S_\Q(\tilde\pi\otimes\eta)$. In addition, $\delta_\infty=-1$ holds
by \eqref{neq1} combined with $\tilde\mu_\infty\eta_\infty=1$ and
$\pi_\infty\cong\pi_\infty^c$.

\section{$4$-dimensional Galois representations of $\Q$}

In the previous section we constructed, for each quadratic character
$\mu$ of $K$ in some dense set $\MM$, a cuspidal automorphic
representation $\Pi^\mu$ of $\GSp_4(\A_\Q)$ such that
$\Pi_\infty^\mu$ is a holomorphic limit of discrete series
representation of weight $(k,2)$. It has the property
\[L^{S(\mu)}(s,\Pi^\mu)=L^{S(\mu)}(s,I_K^\Q(\pi\otimes\mu)),\]
where $S(\mu)$ abbreviates $S_\Q(\pi\otimes\mu)$, $L^{S(\mu)}$
denotes the product of local L-factors outside $S(\mu)$, and
$I_K^\Q$ stands for automorphic induction. Note that $S(\mu)$
includes all the rational primes where $\Pi^\mu$ is ramified.

In this section we shall construct, for each $\mu\in\MM$, a
continuous semisimple representation
\[\tau^\mu: \Gal(\ov\Q/\Q) \to \GL_4(\ov \Q_\ell)\]
such that
\[L^{S(\mu)}(s,\tau^\mu)=L^{S(\mu)}(s,\Pi^\mu).\]
In other words, we shall show that $\pi\otimes\mu$ is associated to
a Galois representation over $\Q$. We shall rely on the following
deep result of Weissauer \cite{We} (see also the closely related
work of Laumon \cite{L1,L3}):

\begin{thm}\label{thm2}
Let $\Pi$ be an irreducible cuspidal automorphic representation of
$\GSp_4(\A_\Q)$ such that $\Pi_\infty$ belongs to the holomorphic
discrete series of weight $(k_1,k_2)$ with $k_1\geq k_2\geq 3$. Let
$S$ denote the union of $\{\ell\}$ and the set of rational primes
where $\Pi$ is ramified. There exists a continuous semisimple
representation
\[\tau: \Gal(\ov\Q/\Q) \to \GL_4(\ov \Q_\ell)\]
such that if $v$ is a rational prime outside $S$ then $\tau$ is
unramified at $v$ and the characteristic polynomial of
$\tau(\Frob_v)$ agrees with the Hecke polynomial of $\Pi$ at $v$. In
other words,
\[L^S(s,\tau)=L^S(s,\Pi),\]
where $L^S$ denotes the product of local $L$-factors outside $S$.
\end{thm}

This theorem was not available for Taylor in \cite{T}; instead, he
utilized the weaker yet powerful results of \cite{T3} to conclude
$L^S(s,\tau)=L^S(s,\Pi)$ for some exceptional set $S$ of zero
Dirichlet density. While the above theorem is not directly
applicable to the representations $\Pi^\mu$ we can combine it with
Taylor's method of pseudo-representations \cite{T2} to achieve our
goal. Our situation is analogous to associating $2$-dimensional
Galois representations to elliptic cusp forms of weight $1$ which
was accomplished by Deligne--Serre in the classical paper \cite{DS}
by a technique involving lifting the weight and then applying a
``horizontal'' family of congruences.

Let $f$ be a Hecke eigenform belonging to $\Pi^\mu$: it is a
vector-valued holomorphic semi-regular Siegel cusp form of weight
$(k,2)$ and some level $N$ (i.e.\ $\Pi^\mu$ is ramified exactly at
the primes dividing $N$). Let $\HH_0^N(\Z)$ be the $\Z$-algebra
generated by the Hecke operators corresponding to primes not
dividing $N$ and denote by $\TT_{(k_1,k_2)}(N)$ the image of
$\HH_0^N(\Z)$ in the space of holomorphic Siegel cusp forms of
weight $(k_1,k_2)$ and level $N$ (see \cite[p.~315]{T2} for precise
definitions). It is known that $\TT_{(k_1,k_2)}(N)\otimes\Q$ is a
semisimple $\Q$-algebra. In particular, we have a homomorphism
$\lambda_f:\TT_{(k,2)}(N)\to\OO_f$ such that $T(f)=\lambda_f(T)f$
for all $T\in\HH_0^N(\Z)$, where $\OO_f$ is the ring of integers of
some (minimally chosen) number field $E_f$.

Using the cup product of $f$ with the $\ell^n$-th power of the
``Hasse Invariant" form exhibited by Blasius and
Ramakrishnan\footnote{Ramakrishnan has pointed out to us that there
is a mistake in \cite{BR}, but it does not affect the part we are
using.} \cite[Proposition~3.6]{BR}, Taylor \cite[Proposition~3]{T2}
constructs a ``vertical" family of morphisms\footnote{We assume here
that $N\geq 3$, otherwise we replace $N$ by $N\ell^2$, say.}
\[\lambda_{n,f}:\TT_{(k,2)+m\ell^n(\ell-1,\ell-1)}(N)\to\OO_f/\ell^{n+1}\]
such that
\[\lambda_{n,f}(T)=\lambda_f(T)\mod\ell^{n+1},\qquad T\in\HH_0^N(\Z).\]
Here $m=m(\ell)$ is a positive integer and $n$ is an arbitrary
positive integer. Together with Theorem~\ref{thm2} this allows us to
apply the theory of pseudo-representations, as in Example~1 in
\cite[\S 1.3]{T2}, to piece together the required Galois
representation $\tau^{\mu}$ corresponding to $\Pi^{\mu}$.

\section{$2$-dimensional Galois representations of $K$}

We have exhibited, using previous notation,
continuous semisimple representations
\[\tau^\mu: \Gal(\ov \Q/\Q) \to \GL_4(\ov \Q_\ell),\qquad\mu\in\MM,\]
satisfying
\begin{equation}\label{eq4}
L^{S(\mu)}(s,\tau^\mu)=L^{S(\mu)}(s,I_K^\Q(\pi\otimes\mu)).
\end{equation}
Note that $S(\mu)$ includes all the rational primes where $\Pi^\mu$
is ramified. If $\chi$ denotes the quadratic character of $\Q$
corresponding to $K$, then we have, by the reciprocity formula,
\[L^{S(\mu)}(s,I_K^{\Q}(\pi\otimes\mu)\otimes\chi)=
L^{S(\mu)}(s,I_K^{\Q}(\pi\otimes\mu)),\] since $\chi$ has trivial
restriction to $K$. This implies, by \eqref{eq4} and the Chebotarev
density theorem, that
\begin{equation}\label{eq5}
\tau^{\mu} \otimes \chi \cong \tau^{\mu}.
\end{equation}
This in turn implies the following
\begin{lem} For each $\mu\in\MM$ there is a continuous semisimple representation
\[ \rho^{\mu}: \Gal(\ov K/K) \to \GL_2(\ov\Q_\ell)\]
such that
\[\tau^{\mu} \cong \Ind_K^{\Q}(\rho^{\mu}).\]
\end{lem}

\begin{proof} Suppose first that $\tau^\mu$ is irreducible over $\Q$.
Then $\Hom(\chi,\tau^{\mu} \otimes (\tau^{\mu})^{\vee} )\neq 0$ by
\eqref{eq5}. Since $\tau^{\mu} \otimes
(\tau^{\mu})^{\vee}=\End(\tau^{\mu})$, we see by Schur's Lemma that
$\tau^{\mu}|_K$ is reducible as $\chi|_K$ is trivial. Let $\rho^\mu$
be an irreducible component of $\tau^\mu|_K$ of minimal dimension
(i.e.\ at most $2$). By Frobenius reciprocity
$\Hom(\Ind_K^\Q(\rho^\mu),\tau^\mu)\neq 0$, hence in fact
$\tau^\mu\cong \Ind_K^\Q(\rho^\mu)$ since $\tau^\mu$ is irreducible
of dimension $4$ and $\Ind_K^\Q(\rho^\mu)$ is of dimension at most
$4$.

Suppose now that $\tau^{\mu}$ is reducible over $\Q$. If $\lambda$
(resp. $\beta$) is a $1$-dimensional (resp. $2$-dimensional)
representation occurring in $\tau^{\mu}$, then \eqref{eq5} shows
that $\lambda \chi$ (resp. $\beta\otimes\chi$) also occurs in
$\tau^{\mu}$. Hence there are four cases to consider:
\begin{enumerate}
\item[1.] $\tau^{\mu}\cong\beta\oplus (\beta\otimes\chi)$.
Then $\tau^{\mu}\cong\Ind_K^{\Q}(\beta|_K)$.
\item[2.] $\tau^{\mu}\cong\beta \oplus \gamma$, where both
$\beta$ and $\gamma$ are $\chi$-invariant. Then $\beta
\cong\Ind_K^{\Q}(\kappa)$ and $\gamma \cong \Ind_K^{\Q}(\nu)$ for
some $1$-dimensional $\kappa$ and $\nu$, so that
$\tau^{\mu}\cong\Ind_K^{\Q}(\kappa\oplus \nu)$.
\item[3.] $\tau^{\mu}\cong\beta \oplus \lambda \oplus \lambda
\chi$. Then $\beta \cong \beta \otimes \chi$, so
$\beta\cong\Ind_K^{\Q}(\kappa)$ for some $1$-dimensional $\kappa$
and $\tau^{\mu}=\Ind_K^{\Q}(\kappa\oplus \lambda|_K)$.
\item [4.]$\tau^{\mu}\cong\lambda \oplus \lambda \chi \oplus
\nu \oplus \nu \chi$. Then $\tau^{\mu}\cong\Ind_K^{\Q}(\lambda|_K
\oplus \nu|_K)$.
\end{enumerate}
\end{proof}

\section{Compatibility of twists}

So far we have constructed, for each quadratic character $\mu$ of
$K$ in some dense set $\MM$, a continuous semisimple representation
\[ \rho^{\mu}: \Gal(\ov K/K) \to \GL_2(\ov\Q_\ell)\] such that
\begin{equation}\label{eqn7}
L^{S(\mu)}(s,\Ind_K^{\Q}(\rho^{\mu}))=L^{S(\mu)}(s,I_K^{\Q}(\pi
\otimes\mu)),\end{equation} where $S(\mu)$ abbreviates
$S_\Q(\pi\otimes\mu)$ and both sides involve Euler factors of degree
$4$ over $\Q$. As we want to compare Euler factors over $K$ it is
useful to rewrite the previous equation (using restriction and base
change) as
\begin{equation}\label{eq6}
L^{S(\mu)}(s,\rho^{\mu})L^{S(\mu)}(s,(\rho^{\mu})^c)=L^{S(\mu)}(s,\pi
\otimes \mu)L^{S(\mu)}(s,(\pi \otimes \mu)^c),
\end{equation} where now $S(\mu)$
abbreviates $S_K(\pi\otimes\mu)$ and all $L$-functions involve Euler
factors of degree $2$ over $K$. Note that $S(\mu)$ includes all the
rational primes where $\rho^\mu$ or $(\rho^\mu)^c$ is ramified.

Our aim is to show that the Galois representations $\rho^\mu$ are
globally compatible in the sense that they can be replaced by twists
$\rho\otimes\mu$ of some fixed $\rho$. This will be achieved in
three lemmata. Recall our assumption that we are not in a case
covered by Remark~\ref{known}.

\begin{lem}\label{nlem1} $\rho^\mu|_L$
is irreducible for all $\mu\in\MM$ and all quadratic extensions
$L/K$.
\end{lem}

\begin{proof} Assume that $\rho^\mu|_L$ is reducible for some $\mu\in\MM$ and
some quadratic extension $L/K$. Let $\Psi$ be an irreducible summand
of $\rho^\mu|_L$ and for a prime $v$ of $K$ outside $S(\mu)$ let
$\{\alpha'_v, \beta'_v\}$ denote the Langlands parameters of $\pi
\otimes \mu$. Implicit in the construction of $\rho^\mu$ is the fact
that it has image in $\GL_2(E)$ for a finite extension of $\Q_\ell$.
In particular, $\Psi:\Gal(\ov L/L) \to E^\times$ is a continuous
character. Applying restriction and base change in \eqref{eq6} we
see that if $w$ is a place of $L$ lying above a place $v$ of $K$
outside $S(\mu)\cup\disc(L/K)$ then
$\Psi(\Frob_w)\in\{(\alpha'_v)^{f}, (\beta'_v)^{f},
(\alpha'_{v^c})^{f}, (\beta'_{v^c})^{f}\}$, where $f=(L_w:K_v)$.
Hence in fact $\Psi(\Frob_w)$ is either one of the Langlands
parameters of the base changes $(\pi\otimes\mu)_L$ or
$(\pi\otimes\mu)_L^c$ at $w$. Applying the results of \cite[\S 3]{T}
we conclude that $(\pi\otimes\mu)_L$ is not cuspidal which by
\cite[Lemma~2]{T} means that we are in Case~1 of Remark~\ref{known}.
This contradiction proves the lemma.
\end{proof}

\begin{lem}\label{twist} Let $\mu\in\MM$ and let $\delta$ be a quadratic character of $K$.
\begin{enumerate}
\item[1.] $\rho^\mu \otimes \delta \not \cong \rho^\mu$ for
$\delta$ nontrivial.
\item[2.] $\rho^\mu \otimes \delta \not \cong (\rho^\mu)^c$ in all cases.
\end{enumerate}
\end{lem}

\begin{proof} Assume first that $\delta$ is nontrivial, and denote
by $L/K$ the corresponding quadratic extension.

Assume that $\rho^\mu \otimes \delta \cong \rho^\mu$. Then
$\Hom(\delta,\rho^\mu \otimes (\rho^\mu)^{\vee} )\neq 0$. Since
$\rho^\mu \otimes (\rho^\mu)^{\vee}=\End(\rho^\mu)$, we see by
Schur's Lemma that $\rho^\mu|_L$ is reducible as $\delta|_L$ is
trivial. This is a contradiction to Lemma~\ref{nlem1} and
establishes the first part of the lemma.

Assume that $\rho^\mu \otimes \delta \cong (\rho^\mu)^c$. Then
$(\rho^\mu)^c\otimes\delta^c\cong\rho$, hence in fact
$\rho^\mu\otimes(\delta\delta^c)\cong\rho$. By the first part of the
lemma this forces $\delta\delta^c=1$, hence $\delta=\delta^c$ since
$\delta$ is quadratic. This implies, using \eqref{eq6}, that
\[L^{T}(s,\rho^{\mu}\otimes\delta)L^{T}(s,(\rho^{\mu}\otimes\delta)^c)=L^{T}(s,\pi
\otimes \mu\delta)L^{T}(s,(\pi \otimes \mu\delta)^c),\] where $T$ is
a finite set of primes in $K$ containing $S(\mu)$. Using again the
assumption $\rho^\mu\otimes\delta\cong(\rho^\mu)^c$ we obtain
\[L^{T}(s,(\rho^\mu)^c)L^{T}(s,\rho^\mu)=L^{T}(s,\pi
\otimes \mu\delta)L^{T}(s,(\pi \otimes \mu\delta)^c),\] hence by
\eqref{eq6}, multiplicity one, and base change, we have in fact
\[I_K^\Q(\pi \otimes \mu\delta)  \cong  I_K^\Q(\pi \otimes
\mu).\] Base changing this to $K$ and comparing the cuspidal
representations (in the isobaric sums), we are forced to have
\begin{equation}\label{neq5}\pi \otimes \mu\delta\cong(\pi\otimes
\mu)^c,\end{equation} since $\pi \otimes \mu\delta\cong\pi\otimes
\mu$ falls under Case~1 of Remark~\ref{known}. Here $\delta =
\delta^c$, so there is a quadratic character $\epsilon$ of $\Q$ such
that $\delta = \epsilon|_K$. Regarding $\delta$ and $\epsilon$ as
id\`ele class characters, $\delta$ is the pull-back of $\epsilon$ by
the norm map $N_{K/\Q}$ and so its restriction to the id\`ele
classes of $\Q$ is the trivial character. Regarding $\delta$ as
Galois character this means that its transfer to $\Gal(\ov \Q/\Q)$
is trivial. Now a standard theorem on Galois cohomology of id\`eles
(cf.\ proof of \cite[Lemma~1]{T}) implies that we may write $\delta
= (\nu/\nu^c)$ for some character $\nu$ of $K$. Plugging this in
\eqref{neq5} we obtain
\[\pi \otimes \mu\nu\cong(\pi\otimes \mu\nu)^c,\]
hence we are in Case~2 of Remark~\ref{known}. This contradiction
establishes the second part of the lemma for nontrivial $\delta$.

It remains to prove the lemma for trivial $\delta$, i.e. that
$\rho^\mu \not\cong (\rho^\mu)^c$. However, this is immediate from
\eqref{eq6} and the multiplicity one theorem since we are assuming
that $\pi \otimes \mu\not\cong(\pi \otimes \mu)^c$.
\end{proof}

\begin{lem} There is a continuous semisimple representation
\[ \rho: \Gal(\ov K/K) \to \GL_2(\ov\Q_\ell)\] which is unramified outside $S$
and for all $\mu\in\MM$ satisfies
\[\rho^\mu\oplus(\rho^\mu)^c \cong
(\rho \otimes \mu) \oplus (\rho \otimes \mu)^c.\]
\end{lem}

\begin{proof}
For any quadratic character $\lambda$ write
$\delta_{\lambda}:=\lambda \lambda^c$. Let $\tilde
\rho^{\mu}:=\rho^{\mu} \otimes \mu^{-1}$. Our goal is to show that
either $\tilde \rho^{\mu}$ or $(\tilde \rho^{\mu})^c \delta_{\mu}$
is independent of $\mu\in\MM$. For any prime $v\notin S$ denote by
$\{\alpha_v, \beta_v\}$ the set of inverse roots of the Hecke
polynomial of $\pi$ at $v$; then for any prime $v \notin S(\mu)$ the
set of eigenvalues of $(\rho^{\mu}\oplus(\rho^{\mu})^c)(\Frob_v)
\cdot (\mu^{-1})(v)=(\tilde \rho^{\mu} \oplus (\tilde \rho^{\mu})^c
\delta_{\mu})(\Frob_v)$ is
$$\{\alpha_v,\beta_v, \alpha_{v^c} \delta_{\mu}(v), \beta_{v^c}
\delta_{\mu}(v) \}$$ by (5.2).

Given $\mu,\mu'\in\MM$ the corresponding set of eigenvalues coincide
over the splitting field $F:=K(\delta_{\mu}\delta_{\mu'})$ of degree
at most $2$. By the Chebotarev density theorem and continuity, so do
their characteristic polynomials. Viewing the representations as
$\overline \Q_l[\Gal(\ov K/K)]$-representations, a general theorem
on semisimple modules of algebras over fields of characteristic 0
\cite[Ch.~8, Sec.~12.2, Prop.~3]{Bour} then tells us that their
semisimplifications are isomorphic. But since the $\rho^{\mu}$ are
semisimple, so is the restriction to the normal subgroup
$\Gal(\overline K/F)$ and we get
\begin{equation}\label{nneq1}(\tilde
\rho^{\mu}\oplus(\tilde \rho^{\mu})^c \delta_{\mu})|_F\cong (\tilde
\rho^{\mu'}\oplus(\tilde \rho^{\mu'})^c
\delta_{\mu'})|_F.\end{equation} By Lemma~\ref{nlem1}, $\tilde
\rho^{\mu}|_F$ is irreducible so either we have $\tilde
\rho^{\mu}|_F \cong \tilde \rho^{\mu'}|_F$ or we have $(\tilde
\rho^{\mu})^c \delta_{\mu}|_F \cong \tilde \rho^{\mu'}|_F$. Let us
fix some element $\mu_0\in\MM$, and for each $\mu\in\MM$ such that
the second case holds for $\mu'=\mu_0$ replace $\tilde \rho^{\mu}$
by $(\tilde \rho^{\mu})^c \delta_{\mu}$ (this corresponds to
replacing the original $\rho^\mu$ by $(\rho^\mu)^c$ which is
legitimate). By this change we have achieved that $\tilde
\rho^{\mu}|_F \cong \tilde \rho^{\mu_0}|_F$ for all $\mu\in\MM$,
that is, $\tilde \rho^{\mu} \cong \tilde \rho^{\mu_0}\otimes
\psi_{\mu,\mu_0}$ for some character $\psi_{\mu,\mu_0}$ of
$\Gal(F/K)$. We shall regard $\psi_{\mu,\mu_0}$ as a quadratic
character of $\Gal(\overline K/K)$ trivial on $\Gal(\overline K/F)$.

Note that $\psi_{\mu_0,\mu_0}=1$, so that for general
$\mu,\mu'\in\MM$ the definition
\begin{equation}\label{nneq3}
\psi_{\mu,\mu'}:=\psi_{\mu,\mu_0}\psi_{\mu',\mu_0}\end{equation} is
unambiguous. This character satisfies
\begin{equation}\label{nneq2}\tilde \rho^{\mu} \cong \tilde \rho^{\mu'}\otimes
\psi_{\mu,\mu'},\end{equation} whence Lemma~\ref{twist} tells us
that in \eqref{nneq1} we must have $\tilde \rho^{\mu}|_F \cong
\tilde \rho^{\mu'}|_F$ for $F=K(\delta_{\mu}\delta_{\mu'})$ and
$\psi_{\mu,\mu'}$ must be trivial on $\Gal(\overline K/F)$. It
follows that $\psi_{\mu,\mu'}=1$ or $\psi_{\mu,\mu'}=\delta_{\mu}
\delta_{\mu'}$, since these are the only characters of
$\Gal(\overline K/K)$ trivial on $\Gal(\overline K/F)$. We claim
that either $\psi_{\mu,\mu'}=1$ for all $\mu,\mu'\in\MM$ or
$\psi_{\mu,\mu'}=\delta_{\mu} \delta_{\mu'}$ for all
$\mu,\mu'\in\MM$. Assume the first alternative fails then there are
$\mu_1,\mu_2\in\MM$ with $\psi_{\mu_1,\mu_2}=\delta_{\mu_1}
\delta_{\mu_2}\neq 1$. For arbitrary $\mu\in\MM$ we have
$\psi_{\mu,\mu_1}\psi_{\mu_,\mu_2}=\psi_{\mu_1,\mu_0}\psi_{\mu_2,\mu_0}=\psi_{\mu_1,\mu_2}$
by \eqref{nneq3}, hence $\psi_{\mu,\mu_1}\neq\psi_{\mu,\mu_2}$.
Therefore $\psi_{\mu,\mu_1}=\delta_{\mu} \delta_{\mu_1}$ or
$\psi_{\mu,\mu_2}=\delta_{\mu} \delta_{\mu_2}$. In fact these
equations are equivalent since the two sides have equal products by
\eqref{nneq3}, namely
\[\psi_{\mu,\mu_1}\psi_{\mu,\mu_2}=\psi_{\mu_1,\mu_0}\psi_{\mu_2,\mu_0}=\psi_{\mu_1,\mu_2}=\delta_{\mu_1}
\delta_{\mu_2}.\] We infer that $\psi_{\mu,\mu_1}=\delta_{\mu}
\delta_{\mu_1}$ is valid for all $\mu\in\MM$. This implies, for all
$\mu,\mu'\in\MM$,
\[\psi_{\mu,\mu'}=\psi_{\mu,\mu_0}\psi_{\mu',\mu_0}=\psi_{\mu,\mu_1}\psi_{\mu',\mu_1}=\delta_{\mu} \delta_{\mu'}. \]

If $\psi_{\mu,\mu'}=1$ for all $\mu,\mu'$ then \eqref{nneq2} implies
that $\tilde \rho^{\mu}$ is independent of $\mu$. If
$\psi_{\mu,\mu'}=\delta_{\mu} \delta_{\mu'}$ for all $\mu,\mu'$ then
\eqref{nneq2} implies that $(\tilde \rho^{\mu})^c \delta_{\mu}$ is
independent of $\mu$. In both cases we denote the common value of
these representations by $\rho$ and verify that it satisfies the
required properties.
\end{proof}

\section{End of proof}

We have shown that there is a dense set $\MM$ of quadratic
characters of $K$ (see Definition~\ref{defn1}) and a continuous
irreducible semisimple representation
\[ \rho: \Gal(\ov K/K) \to \GL_2(\ov\Q_\ell)\] unramified outside $S$ such that
for any $\mu\in\MM$ we have
\begin{equation}\label{eq6b}
L^{S(\mu)}(s,\rho\otimes\mu)L^{S(\mu)}(s,(\rho\otimes\mu)^c)=L^{S(\mu)}(s,\pi
\otimes \mu)L^{S(\mu)}(s,(\pi \otimes \mu)^c),
\end{equation}
where $S(\mu)$ abbreviates $S_K(\pi\otimes\mu)$ (see
Definition~\ref{defn2}). Note that $S(\mu)$ is contained in the
union of $S$ and the set of primes in $K$ where $\mu$ or $\mu^c$ is
ramified.

For any prime $v$ of $K$ outside $S$ denote by $\{\alpha_v,
\beta_v\}$ the set of inverse roots of the Hecke polynomial of $\pi$
at $v$ and by $\{\gamma_v, \delta_v\}$ the inverse roots of the
Frobenius polynomial of $\rho$ at $v$. We shall regard these as
multisets (i.e.\ sets with multiplicities). By \eqref{eq6b} for all
$\mu\in\MM$ unramified at $v$ and $v^c$ we have (as multisets)
\[\{ \gamma_v \mu(v),\delta_v \mu(v),\gamma_{v^c} \mu(v^c),
\delta_{v^c} \mu(v^c) \}=\{\alpha_v \mu(v),\beta_v \mu(v),
\alpha_{v^c} \mu(v^c), \beta_{v^c} \mu(v^c) \}.\] We need to show
that (as multisets)
\begin{equation}\label{end0}\{\gamma_v,\delta_v\}=\{\alpha_v,\beta_v\}\qquad\text{and}\qquad
\{\gamma_{v^c,}\delta_{v^c}\}=\{\alpha_{v^c},\beta_{v^c}\}.\end{equation}

If $v$ is inert then the statement is trivial by the existence of
some $\mu\in\MM$ that is unramified at $v=v^c$.

If $v$ is split then we can find $\mu_1,\mu_2\in\MM$ unramified at
$v$ and $v^c$ such that $\mu_1(v)=\mu_1(v^c)$ but
$\mu_2(v)\neq\mu_2(v^c)$. It follows that (as multisets)
\begin{equation}\label{end4}\{\gamma_v ,\delta_v ,\gamma_{v^c},
\delta_{v^c}\}=\{\alpha_v,\beta_v, \alpha_{v^c},
\beta_{v^c}\}\end{equation} and
\begin{equation}\label{end5}\{\gamma_v ,\delta_v
,-\gamma_{v^c}, -\delta_{v^c}\}=\{\alpha_v,\beta_v, -\alpha_{v^c},
-\beta_{v^c}\}.\end{equation} By forming the sums of both multisets
and then adding and subtracting the two resulting equations we
conclude that
\begin{equation}\label{end1}\gamma_v+\delta_v=\alpha_v+\beta_v\qquad\text{and}\qquad
\gamma_{v^c}+\delta_{v^c}=\alpha_{v^c}+\beta_{v^c}.\end{equation} By
forming the reciprocal sums of both multisets and then adding and
subtracting the two resulting equations we conclude that
\begin{equation}\label{end2}\gamma_v^{-1}+\delta_v^{-1}=\alpha_v^{-1}+\beta_v^{-1}\qquad\text{and}\qquad
\gamma_{v^c}^{-1}+\delta_{v^c}^{-1}=\alpha_{v^c}^{-1}+\beta_{v^c}^{-1}.\end{equation}

Let us focus on the left hand sides of \eqref{end1} and
\eqref{end2}. If the left hand side of \eqref{end2} designates a
nonzero common value then we divide by it the left hand side of
\eqref{end1} and obtain $\gamma_v\delta_v=\alpha_v\beta_v$. Together
with the left hand side of \eqref{end1} this yields
$\{\gamma_v,\delta_v\}=\{\alpha_v,\beta_v\}$, hence in fact the
entire \eqref{end0} upon using \eqref{end4} again. In the same way
\eqref{end0} follows if the right hand side of \eqref{end2}
designates a nonzero common value.

We are left with the subtle case when both sides of \eqref{end2}
designate zero as common value. Then \eqref{end4} and \eqref{end5}
simplify to the same multiset equation
\begin{equation}\label{end6}\{\gamma_v ,-\gamma_v ,\gamma_{v^c},
-\gamma_{v^c}\}=\{\alpha_v,-\alpha_v, \alpha_{v^c},
-\alpha_{v^c}\}\end{equation} and \eqref{end0} simplifies to
\begin{equation}\label{end8}\gamma_v^2=\alpha_v^2\qquad\text{and}\qquad\gamma_{v^c}^2=\alpha_{v^c}^2.
\end{equation} We need to deduce \eqref{end8} from \eqref{end6}.
Taking squares in \eqref{end6} and halving multiplicities we see
that it really is equivalent to the multiset equation
\begin{equation}\label{end9}\{\gamma_v^2,\gamma_{v^c}^2\}=\{\alpha_v^2,\alpha_{v^c}^2\}.\end{equation}
Now we observe that in the present situation
\[\alpha_v^2=-\omega(v)=-\omega(v^c)=\alpha_{v^c}^2,\]
whence in fact \eqref{end9} yields
\[\gamma_v^2=\gamma_{v^c}^2=\alpha_v^2=\alpha_{v^c}^2,\]
so that \eqref{end8} holds as needed.

The proof of Theorem~\ref{thm1} is complete.

\acks{The authors thank the organizers of the Arizona Winter School
2001 where this project was started. The project was designed by
Dinakar Ramakrishnan based on his letter to Peter Sarnak of May
2000. The authors express their special gratitude to Professor
Ramakrishnan for his guidance, generosity, and encouragement without
which this paper would not exist. The second author was supported by
European Community grant MEIF-CT-2006-040371 under the Sixth
Framework Programme.}


\begin{thebibliography}{W}

\bibitem{Bour} N. Bourbaki, \emph{Alg\`{e}bre (\'{E}l\'{e}ments de math\'{e}matique, Fasc. XXIII)}, Hermann, Paris, 1958.

\bibitem{BR} D. Blasius, D. Ramakrishnan, \emph{Maass forms and Galois representations}, Galois groups over $\Q$ (Berkeley, CA, 1987), 33--77, Math. Sci. Res. Inst. Publ., 16, Springer, New York, 1989.

\bibitem{C} L. Clozel, \emph{Motifs et formes automorphes: applications du principe de fonctorialit{\'e}}, Automorphic forms, Shimura varieties, and $L$-functions, Vol. I (Ann Arbor, MI, 1988), 77--159, Perspect. Math., 10, Academic Press, Boston, MA, 1990.

\bibitem{D} P. Deligne, \emph{Formes modulaires et repr{\'e}sentations $\ell$-adiques}, S{\'e}m. Bourbaki 1968/1969, exp. 355, Lecture Notes in Math. 179, Springer-Verlag, 1971, pp. 139--172.

\bibitem{DS} P. Deligne, J.-P. Serre, \emph{Formes modulaires de poids $1$}, Ann. Sci. \'Ecole Norm. Sup. \textbf{7} (1974), 507--530.

\bibitem{FH} S. Friedberg, G. Hoffstein, \emph{Nonvanishing theorems for automorphic L-functions on ${\rm GL}_2$}, Ann. of Math. \textbf{142} (1995), 385--423.

\bibitem{HST} M. Harris, D. Soudry, R. Taylor, \emph{$l$-adic representations associated to modular forms over imaginary quadratic fields I. Lifting to ${\rm GSp}_4(\Q)$}, Invent. Math. \textbf{112} (1993), 377--411.

\bibitem{L1} G. Laumon, \emph{Sur la cohomologie \`a supports compacts des vari\'et\'es de Shimura pour ${\rm GSp}(4)\sb {\Q}$}, Compositio Math. \textbf{105} (1997) 267--359.

\bibitem{L3} G. Laumon, \emph{Fonctions z\^etas des vari\'et\'es de {S}iegel de dimension trois}, Ast\'erisque \textbf{302} (2005), 1--66.

\bibitem{S} J.-P. Serre, \emph{Abelian $l$-adic representations and elliptic curves}, Benjamin, New York, 1968.

\bibitem{T2} R. Taylor, \emph{Galois representations associated to Siegel modular forms of low weight}, Duke Math. J. \textbf{63} (1991), 281--332.

\bibitem{T3} R. Taylor, \emph{On the l-adic cohomology of Siegel threefolds}, Invent. Math. \textbf{114} (1993), 289--310.

\bibitem{T} R. Taylor, \emph{l-adic representations associated to modular forms over imaginary quadratic fields. II}, Invent. Math. \textbf{116} (1994), 619--643.

\bibitem{U} E. Urban, \emph{Sur les repr\'esentations {$p$}-adiques associ\'ees aux repr\'esentations cuspidales de {${\rm GSp}\sb {4/{\Q}}$}}, Ast\'erisque \textbf{302} (2005), 151--176.

\bibitem{We} R. Weissauer, \emph{Four dimensional Galois representations}, Ast\'erisque \textbf{302} (2005), 67--150.

\bibitem{W} A. Wiles, \emph{On ordinary $\lambda$-adic representations associated to modular forms}, Invent. Math. \textbf{94} (1988), 529--573.

\end{thebibliography}
\end{document}